\documentclass{amsart}

\usepackage{amsmath}
\usepackage[dvips]{graphicx}
\usepackage{amsfonts}
\usepackage{amssymb}
\usepackage{amsthm}
\usepackage{mathrsfs}
\usepackage[initials]{amsrefs}
\usepackage[draft]{optional}
\usepackage{graphicx, color}
\usepackage{amssymb, paralist}
\usepackage{amsfonts, amsmath, wasysym}
\usepackage{pdfsync}
\usepackage{latexsym}
\usepackage{graphicx, color}
\usepackage{indentfirst}
\usepackage{bm}
\usepackage[active]{srcltx}

\newcommand{\R}{{\mathbb{R}}}
\theoremstyle{definition}
\newtheorem*{defn}{Definition}
\theoremstyle{definition}
\newtheorem{thm}{Theorem}
\newtheorem{lemma}[thm]{Lemma}
\newtheorem{cor}[thm]{Corollary}

\begin{document}

\title[a conformal positive mass theorem with noncompact boundary]{A CONFORMAL POSITIVE MASS THEOREM WITH NONCOMPACT BOUNDARY}

\author{Alex Freire}
\author{Mohammad Tariquel Islam}

\begin{abstract}
We obtain an integral inequality for asymptotically linear harmonic functions on asymptotically flat 3-manifolds with noncompact boundary, which implies positivity of a convex combination of ADM masses of two conformally related metrics under a positivity condition on a corresponding convex combination of their scalar curvatures. This generalizes a result of Batista and Lopes de Lima, under conditions that do not assume positivity of scalar curvature.
\end{abstract}

\address{Department of Mathematics, University of Tennessee, Knoxville, USA}
\email{afreire@utk.edu}
\email{mislam33@vols.utk.edu}
\maketitle

\section{Introduction}

Our starting point is a result of W.Simon (1999, \cite{simon}) concerning two conformally related, asymptotically flat 3-manifolds with finite ADM mass, $(M,g), (M,g')$, $g'=f^4g$ where $f>0$ is smooth on $M$ and $f$ tends to one at infinity. Using a harmonic spinor argument, W. Simon proves that the ADM masses satisfy $m(g)+m(g')\geq 0$, provided the condition $R+f^4R'\geq 0$ holds pointwise on $M$. Additionally,  $m(g)+m(g')=0$ implies $f\equiv 1$, $M=\R^3$ and $g=g'$ is the euclidean metric. The motivation for this result was an earlier harmonic spinors argument by ul-Alam (1993, \cite{ul-alam}) in the context of a uniqueness result for certain solutions of Einstein's equation coupled with scalar and electric fields, where a Lichnerowicz formula is used without assuming positive scalar curvature of the underlying Riemannian metric. (In particular, ul-Alam obtains the existence of asymptotically constant harmonic spinors  under conditions different from non-negativity of scalar curvature.) 

Almost twenty years later, L-F Tam and Q. Wang (2017, \cite{tam-wang}), in the process of extending W. Simon's result to initial data sets $(M^3, g, k)$ and to asymptotically hyperbolic manifolds, observed that a more general result is a direct consequence of the classical positive mass theorem itself. Namely, for two asymptotically flat, conformally related metrics $g$ and $g'=f^4g$, under the assumption:
$$(1-\lambda)R+\lambda f^4R'\geq 0 \text{ for some }\lambda\in [0,1],$$
one has for the ADM masses: 
$$(1-\lambda)m(g)+\lambda m(g')\geq0;$$
equality implies $M=\R^3$, $f\equiv 1$ and $g=g'$ is the euclidean metric.\vspace{.2cm}

Recall a smooth connected 3-dimensional Riemannnian manifold $(M,g)$ is {\it asymptotically flat} if there exists a compact set $K\subset M$ such that $M\setminus K=\sqcup_{k=1}^NM_{end}^k$, where the ends $M_{end}^k$ are pairwise disjoint and diffeomorphic to the complement of a ball $\R^3\setminus B_1$, and in each end the diffeomorphism defines a coordinate system satisfying:
$$|\partial^l(g_{ij}-\delta_{ij})(x)|=O(|x|^{-q-l}),\text{ for some }\frac 12<q<1\text{ and }l=0,1,2.$$
If we assume the scalar curvature $R_g\in L^1(M,g)$, the ADM mass of each end is well-defined and given by (with $n^0$ the euclidean exterior unit normal to $S_{\rho}$):
$$m(g)=\lim_{\rho\rightarrow \infty}\frac 1{16\pi}\int_{S_{\rho}}\mu_g[n^0]d\sigma^0_{\rho}, \quad \mu_g=\sum_{i,j}(g_{ij,i}-g_{ii,j})dx^j\in \Omega^1(M_{end})$$
To each end of $M$ corresponds an {\it exterior region} $M_{ext}\supset M_{end}$, diffeomorphic to the complement of finitely many balls in $\R^3$ with disjoint closure, with boundary consisting of minimal 2-spheres. \vspace{.2cm}

In Bray et al. (2022, \cite{bkks}), the following is proved: let $(M_{ext},g)$ be an exterior region of a complete asymptotically flat Riemannian 3-manifold with ADM mass $m$. If $u$ is a harmonic function on $(M_{ext},g)$ with Neumann boundary conditions on $\partial M_{ext}$, asymptotic at infinity to one of the coordinate functions in the end considered, we have:
\begin{equation}\label{bkks}
m\geq\frac 1{16\pi}\int_{M_{ext}}\big (\frac{|\nabla^2_gu|^2_g}{|\nabla^g u|_g}+R_g|\nabla^gu|_g\big)dV_g.
\end{equation}
It is noteworthy that the existence of such harmonic functions, with Neumann boundary conditions and asymptotic to a given linear function of the asymptotically flat coordinates, can be proven independently of sign conditions on the scalar curvature. In particular, this inequality proves a strong version of the Positive Mass theorem, and the rigidity statement follows from a short argument (see \cite{bkks}).

It is natural to seek a similar inequality for a convex combination of the ADM masses of conformally related metrics, which would imply the Tam-Wang result alluded to above; we do this in section 2 of this paper. (One issue to keep in mind is that the minimality of $\partial M_{ext}$ with respect to the metric $g$ is used in the proof.)\vspace{.2cm}

The third strand of ideas motivating this work are positive-mass results (and future-timelike energy-momentum vector, for the case of initial data sets) for asymptotically flat manifolds with a noncompact boundary, (\cite{almaraz-barbosa-lima}, \cite{batista-lima}.) We recall the basic definition (we assume for simplicity $M$ has only one end, and $\Sigma$ is connected):

\begin{defn}
A three-dimensional Riemannian manifold $(M,g)$ with non-compact boundary $\Sigma$ is an {\it asymptotically flat half-space} if there exists a compact set $K\subset M$ ($M_{end}=M\setminus K$ and a diffeomorphism $\Phi: M_{end}\rightarrow R^3_+\setminus B_1^+(0)$ such that in the asymptotic coordinate system $y=(y_1,y_2,y_3)$ (with $\R^3_+$ given by $y_3\geq 0$) defined by $\Phi$, the metric $g$ satisfies:
$$|\partial^l(g_{ij}-\delta_{ij})(y)|=O(|y|^{-q-l}),\text{ for some }\frac 12<q<1\text{ and }l=0,1,2.$$ 
We also assume the scalar curvature $R_g$ and the mean curvature $H_g$ of the embedding $\Sigma\hookrightarrow M$ are in $L^1(M,g)$. 
\end{defn}

In this situation, it is proved in \cite{koerber} (Lemma 2.3), assuming also $R_g\geq 0$ and $H_g\geq 0$, that there exists a 3-manifold with boundary $M_{ext}$, diffeomorphic to the complement of finitely many closed balls in $\R^3_+$ and containing an open subset isometric to $M_{end}$, with boundary of the form:
$$\partial M_{ext}=\Sigma_{ext}\sqcup {\mathcal{S}}^{cl}\sqcup {\mathcal{D}}^{fb}.$$
Here $\Sigma_{ext}=\Sigma\cap M_{ext}$, ${\mathcal S}^{cl}$ is a union of disjoint minimal 2-spheres in the interior of $M_{ext}$ and ${\mathcal D}^{fb}$ is a disjoint union of free-boundary minimal 2-disks (i.e. intersecting $\Sigma$ orthogonally.)

An asymptotic invariant of ADM mass type for asymptotically flat half-spaces was defined in \cite{almaraz-barbosa-lima}:
\begin{equation}\label{half-mass-definition}
m_{\Sigma}(g)=\lim_{\rho\rightarrow \infty}\frac 1{16\pi}(\int_{S^+_{\rho}}\mu_g[\mu^0]d\sigma^0_{\rho}+\int_{\partial S^+_{\rho}} \langle \partial_{y_3},\theta^0\rangle_gds^0_{\rho}),
\end{equation}
where $\mu^0$ is the euclidean outer unit normal to the $y$-coordinate hemisphere $S^+_{\rho}$ and $\theta^0$ the outer unit co-normal to the circle $\partial S^+_{\rho}$ in $\Sigma_{ext}$. 

For an asymptotically flat half-space, it is proved in \cite{almaraz-barbosa-lima} (Prop 3.8) that solutions to the following boundary/asymptotic value problem exist:
\begin{equation}\label{harmonic-half-space}
\Delta_gw=0\text{ on }M_{ext};\quad \partial_nw=0\text{ on }{\mathcal S}^{cl}\sqcup {\mathcal D}^{fb},
\end{equation}
$$ w=0\text{ and }|\nabla w|_g\neq 0\text{ on }\Sigma_{ext}; \quad w-y_3=O(|y|^{-q+1+\epsilon})\text{ as}|y|\rightarrow\infty.$$
(This result does not require a sign condition on the scalar curvature.)  Here $n$ is the outward unit normal to $\Sigma_{ext}$, extended to ${\mathcal S}^{cl}\sqcup {\mathcal D}^{fb}$ in the natural way (pointing out of $M_{ext}$). For such harmonic functions $w$, the following inequality is proved in \cite{batista-lima}:
\begin{equation}\label{batista-lima-inequality}
m_{\Sigma}(g)\geq\frac 1{16\pi}\int_{M_{ext}}\big (\frac{|\nabla^2_gw|^2_g}{|\nabla^g w|_g}+R_g|\nabla^gw|_g\big)dV_g+\frac 1{8\pi}\int_{\Sigma_{ext}}H_g|\nabla^gw|_gdA_g.
\end{equation}
Our main observation, proved in section 3, is that this estimate extends to convex combination of masses of conformally related metrics on AF half-spaces.\vspace{.2cm}

\begin{thm}
Let $(M,g)$ be an asymptotically flat 3-dimensional half-space. Suppose $g'=f^4g$ is a conformally related metric on $M$, where $f>0$ and $f\rightarrow 1$ at infinity, so that $g'$ is also asymptotically flat. Assume also $\partial_nf=0$ on $\partial M_{ext}$. Let $\lambda\in [0,1]$. The metric $g_{\lambda}=f^{4\lambda}g$ is also asymptotically flat. Let $w$ be $g_{\lambda}$-harmonic $(\Delta_{\lambda}u=0)$ and satisfy the boundary/asymptotic conditions \eqref{harmonic-half-space}. Then the following inequality holds:
$$
(1-\lambda)m_{\Sigma}(g)+\lambda m_{\Sigma}(g')\geq \frac 1{16\pi}\int_{M_{ext}}\big (\frac{|\nabla^2_{\lambda}w|^2_{\lambda}}{|\nabla^{\lambda} w|_{\lambda}}+((1-\lambda)R_g+\lambda f^4R')|\nabla^{\lambda}w|_{\lambda}\big)dV_{\lambda}$$
$$+\frac 1{8\pi}\int_{\Sigma_{ext}}((1-\lambda)H_g+\lambda f^2H') |\nabla^\lambda w|_{\lambda}dA_{\lambda}.$$
(Subscripts $\lambda$ denote quantities computed in the metric $g_{\lambda}$.)
\end{thm}

\begin{cor}
If $(1-\lambda)R_g+f^4\lambda R'\geq 0$ on $M_{ext}$ and $(1-\lambda)H_g+\lambda f^2 H'\geq 0$ on $\Sigma_{ext}$, we have $(1-\lambda)m_{\Sigma}(g)+\lambda m_{\Sigma}(g')\geq 0$; equality implies $M=\R^3$, $f\equiv 1$ and $g=g'$ is euclidean.
\end{cor}

\vspace{.3cm}

\section{Simon's conformal positive mass theorem via harmonic functions} In this section we describe how the inequality in \cite{bkks} (Theorem 1.2) may be extended to a one-parameter conformal class of asymptotically flat metrics, following an observation in \cite{tam-wang}. We begin by recalling the definition of the {\it exterior region} $M_{ext}$ associated to an end of an asymptotically flat (AF) 3-manifold $(M^3,g)$. 

According to \cite{huisken-ilmanen} (Lemma 4.1), given an AF $(M^3,g)$, associated to each end $M_{end}$ of $M$ there exists an `exterior region' $(M_{ext}, \bar{g})$: a complete, connected, asymptotically flat 3-manifold, with compact boundary $\partial M_{ext}$ consisting of finitely many embedded $g$-minimal ($H_g=0$) 2-spheres, and containing no other minimal surfaces, embedded or immersed. $M_{ext}$ is diffeomorphic to the complement in $\R^3$ of finitely many open balls with disjoint closures. In addition, $(M_{ext},\bar{g})$ contains an open subset isometric to $(M_{end},g)$, and in particular has ADM mass $m(g)$. For this reason, we'll identify $\bar{g}$ and $g$ in the notation from this point on, and regard $M_{end}$ as an open subset of $M_{ext}$, though strictly speaking it isn't. As the authors of \cite{huisken-ilmanen} observe, the existence of $M_{ext}$ is independent of any curvature hypotheses for $(M,g)$. 

We consider a one-parameter family of metrics on $M_{ext}$ conformal to $g$, $g_{\lambda}=f^{4\lambda}g, \lambda\in [0,1]$. The hypotheses on $f$ are such that all $g_{\lambda}$ are AF, namely, as $r=|x|\rightarrow \infty$ on $M_{end}$ (with $x\in \R^3\setminus B_1$ the AF coordinates): 
\begin{equation}\label{conformal-factor}
|f-1|+r|\nabla^gf|_g=O(r^{-q}), \frac 12<q<1; \Delta_gf\in L^1(M,g),\partial_nf=0\text{ on }\partial M_{ext}.
\end{equation}
We need to make sure $\partial M_{ext}$ is minimal with respect to any metric in this conformal family. Recall the transformation formula (where $H_{\lambda}$ is the mean curvature of $\partial M_{ext}$ with respect to the metric $g_{\lambda}$):
$$f^{2\lambda} H_{\lambda}=H_g+4\lambda \frac{\partial_nf}{f},$$
where $n$ is the unit outward normal to $\partial M_{ext}$ (with respect to $g$, or equivalently to any $g_{\lambda}$). Thus the condition $\partial_nf=0$ guarantees $H_{\lambda}=0$. This is the only hypothesis on $\partial M_{ext}$ needed for Thm 1.2 in \cite{bkks}.

The following observation is a special case of one in \cite{tam-wang}

\begin{lemma}\label{mass-g-lambda}
With $f>0$ on $M_{ext}$ as above and $\lambda\in [0,1]$, let $g'=f^4g, g_{\lambda}=f^{4\lambda}g$ (all AF metrics on $M_{ext}$). Then we have for the ADM masses:
$$m(g_{\lambda})=(1-\lambda)m(g)+\lambda m(g').$$
\end{lemma}

\begin{proof}
Recall the expression for the change in ADM mass under conformal change of metric (see e.g. \cite{lee-book}, p. 73):
$$m(g')=m(g)-\frac 1{2\pi}\lim_{\rho\rightarrow \infty}\int_{S_{\rho}}\partial_rfd\sigma^0_{\rho},$$
in particular:
$$m(g_{\lambda})=m(g)-\frac 1{2\pi}\lim_{\rho\rightarrow \infty}\int_{S_{\rho}}\partial_r(f^{\lambda})d\sigma^0_{\rho}=m(g)-\frac {\lambda}{2\pi}\lim_{\rho\rightarrow \infty}\int_{S_{\rho}}\partial_rfd\sigma^0_{\rho},$$
since $f_{|S_{\rho}}\rightarrow 1$ as $\rho\rightarrow \infty$. Thus (using the above for $\lambda=1$):
$$(1-\lambda)m(g)+\lambda m(g')=m(g)+\lambda (m(g')-m(g))=m(g)-\frac {\lambda}{2\pi}\lim_{\rho\rightarrow \infty}\int_{S_{\rho}}\partial_rfd\sigma^0_{\rho}=m(g_{\lambda})$$
\end{proof}
Briefly: the convex combination of ADM masses of two metrics in this family equals the ADM mass of a metric in the family.

To compute the scalar curvature $R_{\lambda}$ of $g_{\lambda}$, recall the standard transformation formula:
$$f^{4\lambda}R_{\lambda}=R_g-8\frac{\Delta_g(f^{\lambda})}{f^{\lambda}},$$
in particular, for $\lambda=1$:
$$-8\frac{\Delta_gf}f=f^4R'-R_g.$$
A simple calculation yields:
$$\frac{\Delta_g(f^{\lambda})}{f^{\lambda}}=\lambda f^{-1}\Delta_gf-\lambda(1-\lambda)f^{-2}|\nabla^gf|^2_g.$$
Thus, we have:
\begin{equation}\label{scalar-lambda}
R_{\lambda}=f^{-4 \lambda}[R_g-8\lambda f^{-1}\Delta_gf+8\lambda(1-\lambda)f^{-2}|\nabla^gf|^2_g]
\end{equation}
$$=f^{-4\lambda}[(1-\lambda)R_g+\lambda f^4R'+8\lambda(1-\lambda)f^{-2}|\nabla^gf|^2_g].$$

Let $M_{ext}(g)$ be the exterior region for the metric $g$ (as described at the start of this section). By the argument in \cite{bartnik} (Theorem 3.1), given $x=(x^1,x^2,x^3)$ AF coordinates on $M_{end}$, for any $\alpha\in S^2\subset\R^3$ there exists $u_{\lambda}$ on $M_{ext}$, $g_{\lambda}$-harmonic ($\Delta_{\lambda}u_{\lambda}=0$) with Neumann boundary conditions $\partial_nu_{\lambda}=0$ on $\partial M_{ext}$, with the asymptotics on $M_{end}$:
$$|u_{\lambda}(x)-\sum_i\alpha_ix^i|=o(|x|^{1-q}),\quad |\partial_iu_{\lambda}-\alpha_i|=o(|x|^{-q}), i=1,2,3.$$
Since, under the assumption $\partial_nf=0$ on $\partial M_{ext}$, we have that $\partial M_{ext}$ is $g_{\lambda}$-minimal, we may apply the inequality \eqref{bkks} to $g_{\lambda}$ (and Lemma 3) to conclude:
\begin{equation}\label{bkks-for-lambda}
(1-\lambda)m(g)+\lambda m(g')=m(g_{\lambda})\geq \frac 1{16\pi}\int_{M_{ext}}\big (\frac{|\nabla^2_{\lambda}u_{\lambda}|^2_{\lambda}}{|\nabla^{\lambda} u_{\lambda}|_{\lambda}}+R_{\lambda}|\nabla^{\lambda}u_{\lambda}|_{\lambda}\big)dV_{\lambda},
\end{equation}
where $R_{\lambda}$ is given by \eqref{scalar-lambda}.

A short argument then leads to a proof of the main result in \cite{simon}:
\begin{thm}
Let $(M^3,g)$ be asymptotically flat. Let $f>0$ on $M_{ext}(g)$ satisfy \eqref{conformal-factor}, so that $g'=f^4g$ is also AF.  Let $\lambda\in [0,1]$. Assume $(1-\lambda) R_g+\lambda f^4 R'\geq 0$ on $M_{ext}(g)$. Then: (i) $(1-\lambda)m(g)+\lambda m(g')\geq 0$. (ii)  If $(1-\lambda)m(g)+\lambda m(g')=0$, then $f\equiv 1, M^3=\mathbb{R}^3$ and $g'=g$ is euclidean.
\end{thm}

\begin{proof}
(i) Follows directly from \eqref{bkks-for-lambda}; note $R_{\lambda}\geq 0$, from expression \eqref{scalar-lambda} above.\vspace{.2cm}

(ii) Suppose $(1-\lambda)m+\lambda m(g')=0$. Then $m(g_{\lambda})=0$, so by uniqueness in the positive mass theorem we have $M=\mathbb{R}^3$ and $g_{\lambda}$ is euclidean. From expression \eqref{scalar-lambda} and the assumption $(1-\lambda)R_g+\lambda f^4R'\geq 0$, it follows that $\nabla_gf\equiv 0$ on $M$, and since $f\rightarrow 1$ at infinity this implies $f\equiv 1$. So $g=g'=g_{\lambda}$, and $g$ is euclidean.
\end{proof}

\section{Proof of Theorem 1 and Corollary 2.}

\begin{proof} First observe that, on the boundary $\partial S_{\rho}^+$ of the coordinate hemisphere, we have:
$$\langle \partial_{y_3},\theta^0\rangle_{g_{\lambda}}=f^{4\lambda} \langle \partial_{y_3},\theta^0\rangle_g\rightarrow \langle \partial_{y_3},\theta^0\rangle_g,$$
as $\rho\rightarrow \infty$. \vspace{.2cm}
Thus the proof of Lemma 3 still applies in the half-space case: there is no contribution from $\partial S_{\rho}^+$ in the equality:
$$m_{\Sigma}(g_{\lambda})=m_{\Sigma}(g)-\frac{\lambda}{2\pi}\lim_{\rho}\int_{S_{\rho}^+}\partial_rfd\sigma_{\rho}^0.$$
This implies:
$$(1-\lambda)m_{\Sigma}(g)+\lambda m_{\Sigma}(g')=m_{\Sigma}(g)+\lambda(m_{\Sigma}(g')-m_{\Sigma}(g))$$
$$=m_{\Sigma}(g)-\frac{\lambda}{2\pi} \lim_{\rho}\int_{S_{\rho}^+}\partial_rfd\sigma^0_{\rho}=m_{\Sigma}(g_{\lambda}),$$
exactly as in Lemma 3.

\noindent From the mean curvature transformation formula:
\begin{equation} H_{\lambda} = f^{-2\lambda}H_g +4 \lambda f^{-2\lambda-1}\partial_nf,\end{equation}
for $\lambda=1$ we get 
\[ H'= f^{-2}H_g + 4 f^{-3}\partial_nf,\]or
\begin{equation}\partial_nf = \frac{1}{4}\left(f^3H'-fH_g\right). \end{equation}
Now replacing $\partial_nf$ in equation (8) by (9) we find:
\begin{equation}
	H_{\lambda}=f^{-2\lambda}\left[(1-\lambda)H_g+\lambda f^2H'\right].
\end{equation}

\noindent Also the scalar curvature transformation formula is \begin{equation}
	R_{\lambda} = f^{-4\lambda}\left[(1-\lambda)R_g + \lambda f^4 R' + 8\lambda(1-\lambda)f^{-2}\left|\nabla^g f\right|^2_g\right].
\end{equation}

\noindent Thus for $m_{\Sigma}(g_{\lambda})$  the inequality in \cite{batista-lima} (Theorem 1.6) becomes:

\[\begin{aligned}
	(1-\lambda)m_{\Sigma}(g)+\lambda m_{\Sigma}(g')= m_{\Sigma}(g_{\lambda}) &\geq \frac{1}{16\pi}\int_{M_{\text{ext}}}\left(\frac{\left|\nabla^2_{\lambda} u\right|^2_{\lambda}}{\left|\nabla^{\lambda} u\right|_{\lambda}}+R_{\lambda}\left|\nabla^{\lambda} u\right|_{\lambda}\right)dV_{\lambda} + \frac{1}{8\pi}\int_{\Sigma_{\text{ext}}} H_{\lambda}\left|\nabla^{\lambda} u\right|_{\lambda} dA_{\lambda}\\
		&\hspace{-1.5in}= \frac{1}{16\pi}\int_{M_{\text{ext}}}\left(\frac{\left|\nabla^2_{\lambda} u\right|^2_{\lambda}}{\left|\nabla^{\lambda} u\right|_{\lambda}}+\left[f^{-4\lambda}\left\{(1-\lambda)R_g + \lambda f^4 R' + 8\lambda(1-\lambda)f^{-2}\left|\nabla^gf\right|^2_g  \right\} \right]\left|\nabla^{\lambda} u\right|_{\lambda}\right)dV_{\lambda} \\
		&\hspace{.5in}+ \frac{1}{8\pi}\int_{\Sigma_{\text{ext}}} f^{-2\lambda}\left[(1-\lambda)H_g+\lambda f^2H'\right]\left|\nabla^{\lambda} u\right|_{\lambda} dA_{\lambda}
\end{aligned}\]

Now the inequality in Theorem 1 follows from noting that all the terms on the right hand side of the equality are non-negative. This inequality immediately implies the first statement on Corollary 2. 
\end{proof}

\noindent Turning to the proof of the rigidity statement:

\begin{proof}
\noindent If $\displaystyle (1-\lambda)m_{\Sigma}(g)+\lambda m_{\Sigma}(g') =0$, it follows from the inequality in the above proof that:
$$(1-\lambda)R_g + \lambda f^4 R'   = 0, \quad (1-\lambda)H_g+\lambda f^2H' = 0 \text{ and }\nabla^2_{\lambda} u \equiv 0,$$
as well as, using (11):
$$\nabla^gf\equiv 0\text{ and }R_{\lambda}\equiv 0.$$
\noindent Since $f\to 1$, $f\equiv 1.$  Thus $g=g'=g_{\lambda}$. In addition, for the second fundamental form $A$ of $\Sigma$ in $M$:
$$A=\frac{\nabla_{\Sigma(\lambda)}^2 u}{\left|\nabla_{\lambda}u\right|}\equiv 0,$$
where $\nabla_{\Sigma(\lambda)}^2 u$ denotes the Hessian $\nabla^2_{\lambda}u$ of $u$ restricted to $T\Sigma\otimes T\Sigma$.  Thus $\Sigma_{\text{ext}}$ is totally geodesic.\vspace{.2cm}

At this point one can simply use the proof of the rigidity part of Theorem 1.4 in \cite{batista-lima} to conclude. For the reader's convenience, we outline the argument. Reflect $\left(M_{\text{ext}},g_{\lambda}\right)$ across the totally geodesic boundary $\Sigma_{ext}$  to obtain its double,  a smooth $3-$manifold $\left(\tilde{M}_{\text{ext}},\tilde{g}_{\lambda}\right)$ which is also asymptotically flat, with boundary consisting of finitely many minimal two-spheres. It follows from the asymptotics of $(\tilde{M}_{ext}, \tilde{g}_{\lambda})$ that:
$$m_{ADM}(\tilde{M}_{ext}, \tilde{g}_{\lambda})= 2 m_{\Sigma}(g_{\lambda})=0.$$
Now the argument in the last paragraph of Section 6.2 in \cite{bkks} (based on uniqueness in the Positive Mass Theorem) implies there are, in fact, no minimal two-spheres, and hence $\left(\tilde{M},\tilde{g}_{\lambda}\right)=\left(\mathbb{R}^3,\delta\right)$, which implies $\left(M,g_{\lambda}\right)=\left(\mathbb{R}^3_+,\delta\right)$.
\end{proof}

\section{Appendix: application to a theorem of ul-Alam}

In \cite{ul-alam}, ul-Alam considers a static Lorentzian metric on $M_{ext}\times \R$, of the form $-V^2dt^2+g$. Here $(M_{ext},g)$ is asymptotically flat with minimal boundary, diffeomorphic to the complement of a closed ball in $\R^3$, and $V$ (the `static potential') is a positive function on $M_{ext}$. The author considers solutions of the Einstein equations for $g$, coupled with a scalar field $\phi$ and an electric potential $U$ (functions on $M_{ext}$). A solution of the system has scalar curvature and static potential satisfying:

\begin{equation}\label{static-potential}
R_g=2V^{-2}e^{-2\phi}|\nabla U|^2+2|\nabla \phi|^2,\quad \Delta_gV=V^{-1}e^{-2\phi}|\nabla U|^2
\end{equation}
The scalar field $\phi$ and electric potential $U$ satisfy:
\begin{equation}\label{scalar-field-electric-potential}
\Delta_gU=V^{-1}\langle \nabla V,\nabla U\rangle_g+2\langle \nabla \phi,\nabla U\rangle_g,\quad \Delta_g\phi=-V^{-1}\langle \nabla V,\nabla \phi\rangle_g+V^{-2}e^{-2\phi}|\nabla U|_g^2
\end{equation}
The following asymptotics are assumed as $r=|x|\rightarrow \infty$ in AF coordinates:
\begin{equation}\label{asymptotics}
g=(1+\frac{2M}r)\delta+O_2(r^{-2}), \quad V=1-\frac Mr+O_2(r^{-2}),
\end{equation}
$$e^\phi=1-\frac{Q^2}{2Mr}+O_2(r^{-2}), \quad U=\frac Qr+O_2(r^{-2}).$$
Here $M>0$ and $Q$ are constants, with $Q^2<2M^2$. In particular, $m_{ADM}(g)=M.$
At $\partial M_{ext}$ (corresponding to $r=2M$), the following boundary conditions are assumed:
$$\langle \nabla V,\nabla \phi\rangle_g=0,\quad \langle \nabla V,\nabla U\rangle_g=0,\quad |\nabla V|_g^2=const.\quad V=0.$$
The main result of \cite{ul-alam} is that the unique solution to the system (including the asymptotics and the boundary conditions) is the one found by Gibbons:
$$g=V^{-2}dr^2+(1-\frac{Q^2}{Mr})r^2d\omega^2, \quad V^2=1-\frac{2M}r,\quad e^{2\phi}=1-\frac {Q^2}{Mr},\quad U=\frac Qr.$$
The strategy of proof is to use conformal deformations of $g$ and appeal to an extension of the harmonic spinor proof of the Positive Mass Theorem (in the sense that nonnegative scalar curvature is not assumed; a proof of existence of harmonic spinors without this assumption is given in \cite{ul-alam}.) Here we describe how the argument can be seen as an application of Theorem 4 in section 2.\vspace{.2cm}

Two conformal factors are defined on $M_{ext}$:
$$\chi=\frac 12e^{-\phi/2}[(1+Ve^\phi)^2-2U^2]^{1/2},\quad \psi=\frac 12 e^{-\phi/2}(e^{\phi}+V).$$
It is convenient to consider the functions on $M_{ext}$:
$$u=1+Ve^{\phi}+\sqrt{2}U,\quad v=1+Ve^{\phi}-\sqrt{2}U,$$
so we have: $\chi=(1/2)e^{-\phi/2}\sqrt{uv}$.

From the equations for $V,\phi$ and $U$ one obtains (Lemma 1 in \cite{ul-alam}):
\begin{equation}\label{u-v-equations}
\Delta_gu-\langle \nabla u, X-\nabla \phi\rangle_g=0,\quad \Delta_gv+\langle \nabla v,X-\nabla \phi\rangle_g=0, \text{ where }X=\sqrt{2}V^{-1}e^{-2\phi}\nabla U.
\end{equation}
It is proved in Lemmas 2 and 3 of \cite{ul-alam} that $u>0,v>0,e^{\phi}+V>0$ on $M_{ext}$. Thus, $\chi>0, \psi>0$ on $M_{ext}$ can be used as conformal factors. Define:
$$\gamma=\chi^4g, \quad \eta=\psi^4g.$$
{\it Remark:} In \cite{ul-alam}, a reflection argument is used (across the minimal boundary of $M_{ext}$), and the functions $\chi, \psi$ are extended to the reflected 3-manifold $M_{ext}^-$, in such a way that $\chi^4g$ `compactifies' to a metric $\gamma$ on a 3-manifold without boundary $(N,\gamma)$, complete and AF at infinity on the original $M_{ext}$. The goal is to show $(N,\gamma)$ is isometric to $(\R^3, \delta)$; the positive mass theorem cannot be invoked directly, since $R_{\gamma}\geq 0$ is not clear. Here we concentrate on $M_{ext}$ to reduce this result to the uniqueness statement in Theorem 1.\vspace{.2cm}

{\bf Lemma 1.}\footnote{ This is Lemma 4 in \cite{ul-alam}, but a detailed proof is not included in that paper.}
$$(i) \psi^4R_{\eta}=e^{2\phi}R_{e^{2\phi}g},\quad (ii)\chi^4R_{\gamma}=2\mathcal{P}-e^{2\phi}R_{e^{2\phi}g}, $$
where $\mathcal{P}=|\frac{\nabla u}u-\frac{\nabla v}v-X|^2$.
\begin{proof}
(i) From the expression for transformation of scalar curvature under conformal change, claim (i) is equivalent to:
$$\frac{\Delta_g\psi}{\psi}=\frac{\Delta_g(e^{\phi/2})}{e^{\phi/2}}=\frac{\Delta_g\phi}2+\frac{|\nabla \phi|^2_g}4.$$
By direct calculation from the definition of $\psi$:
$$2\Delta_g\psi=\Delta_g(e^{\phi/2})+Ve^{\phi/2}[V^{-1}\Delta_gV-V^{-1}\langle \nabla \phi,\nabla V\rangle_g-\frac{\Delta_g\phi}2+\frac{|\nabla \phi|^2_g}4].$$
From the field equations \eqref{static-potential},\eqref{scalar-field-electric-potential} for $\phi$ and $V$, we have:
$$\Delta_g\phi=V^{-1}\Delta_gV-V^{-1}\langle \nabla \phi,\nabla V\rangle_g.$$
Thus:
$$2\Delta_g\psi=(e^{\phi/2}+Ve^{\phi/2})(\frac{\Delta_g}2+\frac{|\nabla \phi|^2_g}4)=2\psi(\frac{\Delta_g\phi}2+\frac{|\nabla \phi|^2_g}4),$$
as we wished to show.

(ii) Let $\chi=e^w$, where $w=-\log 2+\frac 12(\log u+\log v-\phi)$, so $\chi^{-1}\Delta_g\chi=\Delta_gw+|\nabla w|^2_g$. Claim (ii) is equivalent to:
$$R_g-8\frac{\Delta_g\chi}{\chi}=2\mathcal{P}-(R_g-8\frac{\Delta_g(e^{\phi/2})}{e^{\phi/2}}),$$
or to (using the expression for $R_g$ from \eqref{static-potential}):
$$|X|^2_g+2|\nabla \phi|_g^2=4(\Delta_gw+|\nabla w|_g^2+\frac{\Delta_g\phi}2+\frac{|\nabla \phi|^2_g}4)+\mathcal{P}.$$
This follows by direct calculation, using \eqref{u-v-equations}.
\end{proof}
Consider the function $f=\frac{\chi}{\psi}$, so that $\gamma=f^4\eta$. Combining conclusions (i) and (ii) of the lemma, we have:
$$\psi^4R_{\eta}+\chi^4R_{\gamma}=2\mathcal P\geq 0,\text{ or }R_{\eta}+f^4R_{\gamma}=2\mathcal{P}\psi^{-4}\geq 0.$$
As observed in \cite{ul-alam} (Lemma 8), the asymptotics in \eqref{asymptotics} imply for $f$:
$$f=1+O(r^{-2}),\quad \nabla f=O(r^{-3})\text{ as }r\rightarrow \infty.$$
For the ADM masses of $\gamma$ and $\eta$, it is easy to see we have:
$$m_{ADM}(\gamma)=m_{ADM}(\eta)=0.$$
For $\gamma$ this follows from the well-known expression:
$$m_{ADM}(\chi^4g)=m_{ADM}(g)-\frac 1{2\pi}\int_{S_{\infty}}\partial_r\chi d\sigma^0=m_{ADM}(g)-M=0,$$
from the asymptotics of $\psi$ that follow from \eqref{asymptotics}, namely:
$$\psi=1-\frac{M}{2r}+O_1(r^{-2}),\quad \partial_r\psi=\frac{M}{2r^2}+O(r^{-3}).$$
And then the conclusion for $\eta$ follows from the same expression and the asymptotics of $f$, since $\eta=f^{-4}\gamma$. 

Thus we are in the situation of uniqueness in Theorem 4, in the special case where $m(g)=m(g')=0$. This leads to a short, spinor-free proof of the following result in \cite{ul-alam}, the main ingredient in the uniqueness theorem that is the subject of that paper.

{\bf Theorem.} (\cite{ul-alam})

 We have $\mathcal{P}\equiv 0, f\equiv 1$, $\gamma=\eta$, $(N,\gamma)=(\R^3,\delta)$.

{\it Proof.} This follows directly from the uniqueness statement (ii) in Theorem 4 (with $\lambda=1/2$), which applies since $R_{\eta}+f^4R_{\gamma}=2\mathcal{P}\psi^{-4}\geq 0$.

\end{document}